\documentclass[11pt]{article} \textheight = 24truecm \textwidth
= 16truecm \hoffset = -2truecm \voffset = -2truecm
\usepackage[english]{babel}
\usepackage{amssymb,amsthm}
\usepackage{amsmath}

\begin{document}

\large \sloppy

\begin{center}
\textsc {On a metric on the space of idempotent probability
measures}

\textsl{A. A. Zaitov, I. I. Tojiev}
\end{center}

\begin{abstract}
In this paper we construct a metric on the space of idempotent
probability measures on the given compactum, which is an
idempotent analog of the Kantorovich metric on the space of
probability measures.
\end{abstract}



\section{Introduction}

Our goal is to prove the following

$\textbf{Theorem 1.}$ \textsl{The function $\rho_I:I(X)\times I(X)
\rightarrow R$ is a metric on the space $I(X)$ of idempotent
probability measures on the given compactum $X$ which generates
pointwise convergence topology in $I(X)$, where $\rho_I$ defines
by the equality (6).}

A metrizable compact Hausdorff space is called a compactum.

The metric $\rho_I$ on the space $I(X)$ of idempotent probability
measures on the given compactum $X$ constructed in present paper
is an positive answer to the Question 7.8 placed in [1] (the
question expresses as follows: is there a counterpart of the
Prokhorov metric for the functor of idempotent probability
measures?)

The notion of idempotent measure finds important applications in
different part of mathematics, mathematical physics, economics,
mathematical biology and others. One can find a row of
applications of idempotent mathematics from [2].

Let $S$ be a set equipped with two algebraic operation: addition
$\oplus$ and multiplication $\odot$. $S$ is called a semiring if
the following conditions hold:

($i$) the addition $\oplus$ and the multiplication $\odot$ are
associative;

($ii$) the addition $\oplus$ is commutative;

($iii$) the multiplication $\odot$ is distributive with respect to
the addition $\oplus$.

A semiring $S$ is commutative if the multiplication $\odot$ is
commutative. A unity of semiring $S$ is an element $\textbf{1}\in
S$ such that $\textbf{1}\odot x=x\odot\textbf{1}=x$ for all $x\in
S$. A zero of a semiring $S$ is an element $\textbf{0}\in S$ such
that $\textbf {0}\neq \textbf{1}$ and $\textbf{0}\oplus x=x$,
$\textbf{0}\odot x=x\odot \textbf{0}=\textbf{0}$ for all $x\in S$.
A semiring $S$ is idempotent if $x\oplus x=x$ for all $x\in S$. A
semiring $S$ with zero $\textbf{0}$ and unity $\textbf{1}$ is
called a semifield if each nonzero element $x\in S$ is invertible.

Let ($S, \oplus, \odot, \textbf{0},\textbf{1}$) be a semiring. On
$S$ a partially order $\prec$ arises by naturally way: for
elements $a, b\in S$ by definition we have $a\prec b$ if and only
if $a\oplus b=b$. So, all elements of $S$ are nonnegative:
$\textbf{0}\prec x$ for all $x\in S$. The idempotent analog of
functions are maps $X\rightarrow S$ where $X$ is an arbitrary set
and $S$ is an idempotent semiring. $S$-valued functions may be
added, multiplied by each other and multiplied by elements of $S$.

The idempotent analog of the linear space of functions is a set of
$S$-valued functions $X\rightarrow S$, which is closed under
addition of functions and multiplication of functions by elements
of $S$ (which is $S$-semimodule). Denote by $B(X,S)$ the
semimodule of functions $X\rightarrow S$ that are bounded in the
sense of the order $\prec$ on $S$. A functional $f:B(X,
S)\rightarrow S$ by definition is idempotent linear (or
\textit{maxplus}-linear) if
$$
f(\lambda_{1}\odot\varphi_{1}\oplus \lambda_{2}\odot
\varphi_2)=\lambda_{1}\odot f(\varphi_1)\oplus \lambda_{2}\odot
f(\varphi_{2})
$$
for all  $\lambda_{1}$, $\lambda_{2}\in S$ and $\varphi_{1}$,
$\varphi_{2}\in B(X, S)$.

Let $\mathbb{R}$ be the field of real numbers and $\mathbb{R_{+}}$
the semifield of nonnegative real numbers (with respect to the
usual operations). The change of variables $x\mapsto u=h\ln x$,
$h>0$, defines a map $\Phi_{h}:\mathbb{R_{+}}\rightarrow
S=\mathbb{R}\cup\{-\infty\}$. Let the operations of addition
$\oplus$ and multiplication $\odot$ on $S$ be the images of the
usual operations of addition $+$ and multiplication $\cdot$ on
$\mathbb{R}$, respectively, by the map $\Phi_{h}$, i. e. let
$u\oplus_{h}v=h \ln(\exp (u/h)+\exp(v/h)),$ $u\odot v=u+v$. Then
we have $\textbf{0}=-\infty=\Phi_{h}(0)$,
$\textbf{1}=0=\Phi_{h}(1)$. It is easy to see that
$u\oplus_{h}v\rightarrow \max\{u, v\}$ as $h\rightarrow 0$. Hence,
$S$ forms semifield with respect to operations $u\oplus
v=\textrm{max}\{u,v\}$ and $u\odot v=u+v$. It denotes by
$\mathbb{R}_{\textrm{max}}$. It is idempotent. This passage from
$\mathbb{R}_+$ to $\mathbb{R}_{\textrm{max}}$ is called the Maslov
dequantization.

Let $X$ be a compact Hausdorff space, $C(X)$ the algebra of
continuous functions $\varphi : X \rightarrow \mathbb{R}$ with the
usual algebraic operations. On $C(X)$ the operations $\oplus$ and
$\odot$ define as follow:

$\varphi \oplus \psi=\textrm{max}\{\varphi, \psi \}$, where
$\varphi, \psi \in C(X)$,

$\varphi \odot \psi=\varphi + \psi$, where $\varphi$, $\psi \in
C(X)$,

$\lambda \odot \varphi=\varphi+\lambda_{X}$, where $\varphi\in
C(X)$, $\lambda\in \mathbb{R},$  and $\lambda_X$ is a constant
function.

Recall [1] that a functional $\mu : C(X)\rightarrow
\mathbb{R}(\subset \mathbb{R}_{\textrm{max}})$ is called to be an
idempotent probability measure on $X$, if:

1)  $\mu (\lambda_{X})=\lambda$ for each $\lambda \in \mathbb{R}$;

2)  $\mu (\lambda \odot \varphi)=\mu (\varphi)+\lambda$ for all
$\lambda \in \mathbb{R}$, $\varphi \in C(X)$;

3)  $\mu (\varphi \oplus \psi)=\mu(\varphi)\oplus \mu(\psi)$ for
every $\varphi$, $\psi\in C(X)$.

The number $\mu(\varphi)$ is named Maslov integral of the function
$\varphi\in C(X)$ with respect to $\mu$.

For a compact Hausdorff space $X$ a set of all idempotent
probability measures on $X$ denotes by $I(X)$.  Consider $I(X)$ as
a subspace of $\mathbb{R}^{C(X)}$. In the induced topology the
sets
\begin{center} $\langle \mu; \varphi_1, \varphi_2, ...,
\varphi_k; \varepsilon \rangle=\{\nu\in I(X):
|\mu(\varphi_i)-\nu(\varphi_i)|<\varepsilon, i=1, ..., k\}$,
\end{center}
form a base of neighborhoods of the idempotent measure $\mu\in
I(X)$, where $\varphi_i\in C(X)$, $i=1, ..., k$, and $\varepsilon
>0$. The topology generated by this base coincide with pointwise
topology on $I(X)$. The topological space $I(X)$ is compact [1].

Given a map $f:X\rightarrow Y$ of compact Hausdorff spaces the map
$I(f):I(X)\rightarrow I(Y)$ defines by the formula
$I(f)(\mu)(\varphi)=\mu(\varphi\circ f)$, $\mu\in I(X)$, where
$\varphi\in C(Y)$.

The construction $I$ is a covariant functor, acting in the
category of compact Hausdorff spaces and their continuous
mappings. Moreover, $I$ is normal functor [1]. Hence, in
particular, yields that if $X$ is compactum then $I(X)$ is also
compactum. From here a question about construction a metric on
$I(X)$ generated the pointwise topology on $I(X)$ arises in a
naturally way. In the present paper we will construct this metric.

Since $I$ is normal functor then for an arbitrary idempotent
measure $\mu\in I(X)$ we may define the support of $\mu$:
\begin{center}supp $\mu=\bigcap\{A\subset X:
\overline{A}=A,$ $\mu\in I(A)\}$.
\end{center}
For brevity, put $S\mu=\mbox{supp}\mu$.

For a positive integer $n$ put $I_n(X)=\{\mu\in I(X):
|\textrm{supp} \mu|\leq n\}$. Define the following set
\begin{center}
$I_{\omega}(X)=\bigcup\limits_{n=1}^\infty I_n(X)$.
\end{center}
It is known [1] that $I_{\omega}(X)$ is everywhere dense in
$I(X)$. A functional $\mu\in I_{\omega}(X)$ is named as an
idempotent probability measure with finite support.

Note that if $\mu$ is an idempotent probability measure with a
finite support $\textrm{supp} \mu=\{x_1, x_2, ..., x_k \}$ then it
represents in the form
$$
\mu=\lambda_1\odot\delta_{x_1}\oplus\lambda_2\odot\delta_{x_2}\oplus
...\oplus\lambda_k\odot\delta_{x_k} \eqno(2)
$$
uniquely, where $\lambda_i\in \mathbb{R}_{\textrm{max}}$, $i=1,
..., k$, $\lambda_1\oplus\lambda_2\oplus
...\oplus\lambda_k=\textbf{1}$.

\section{The proof of Theorem 1}

Let $\mu_1,\ \mu_2\in I_{\omega}(X)$. Then by (2) we have
$\mu_k=\bigoplus\limits_{i=1}^{n_k}\lambda
_{k}\odot\delta_{x_{ki}}$, $i=1,\ 2$. Put
$$\Lambda_{12}=\Lambda(\mu_1,\ \mu_2)=\{\xi\in I(X^2): I(\pi_i)(\xi)=\mu_i,\ i=1,2\},$$
where $\pi_i:X\times X\longrightarrow X$ is projection onto $i$-th
faction, $i=1,2$. We show the set $\Lambda(\mu_1,\ \mu_2)$ is
nonempty. Without loosing of generality suppose
$\lambda_{11}=\lambda_{21}=0$. Directly checking then shows that
$I(\pi_i)(\xi)=\mu_i$, $i=1,\ 2$, for all idempotent probability
measures $\xi\in I(X^2)$ of the form $\xi=\xi^0\oplus R(\mu_1,\
\mu_2)$. Here
$$\xi^0=0\odot\delta(x_{11},x_{21})\oplus\bigoplus\limits_{t=2}^{n_2}\lambda_{2t}\odot\delta(x_{11},x_{2t})
\oplus\bigoplus\limits_{t=2}^{n_1}\lambda_{1t}\odot\delta(x_{1t},x_{21})$$
is idempotent probability measures on $X^2$ and
$$
R(\mu_1,\mu_2)=\bigoplus_{\footnotesize
\begin{array}{c}
{k\in K},\\
{m\in M}
\end{array}}\gamma_{km}\odot\delta(x_{1k},x_{2m})
$$
is some functional on $C(X)$ where $$\gamma_{km}\leq
\text{min}\{\lambda_{1k}, \lambda_{2m}\},\ k\in K,\ m\in M,\
K\subset\{2,..., n_1\},\ M\subset\{2,..., n_2\}.$$ Thus
$\xi\in\Lambda(\mu_1,\ \mu_2)$. Here in fact more is proved: it is
easy to see if $n_1\geq 2$ and $n_2\geq 2$ then quantity of the
numbers $\gamma_{km}$ is uncountable. From here one concludes that
the cardinality of the set $\Lambda(\mu_1,\ \mu_2)$ is no less
than continuum as soon as each of supports $S\mu_i$, $i$=1,2,
contains no less than two points.

Note that $\xi=\xi^0$ if one takes empty set as $K$ and $M$.

By definition for each idempotent probability measure
$\xi\in\Lambda(\mu_1,\ \mu_2)$ we have $\bigoplus\limits_{(x_j,
x_k)\in S\xi}|\lambda_{2k}-\lambda_{1j}|\odot\rho(x_{1j},
x_{2k})<\infty$. In other hand as the set
$$\{|\lambda_{2k}-\lambda_{1j}|\odot\rho(x_{1j}, x_{2k}): j=1,...,
n_1; k=1,..., n_2\}\eqno{(3)}$$ is finite there exists the number
$$
\min_{\xi\in\Lambda_{12}}\left\{\bigoplus\limits_{(x_j, x_k)\in
S\xi} |\lambda_{2k}-\lambda_{1j}|\odot\rho(x_{1j},
x_{2k})\right\}.
$$

Put
$$H(\mu_1,\ \mu_2)=\min_{\xi\in\Lambda_{12}}\left\{\bigoplus\limits_{(x_j, x_k)\in S\xi}
|\lambda_{2k}-\lambda_{1j}|\odot\rho(x_{1j},
x_{2k})\right\}.\eqno{(4)}$$

Here the following statement takes place.

\textbf{Lemma 1.} \textsl{For an arbitrary pair $\mu_1,\ \mu_2\in
I_{\omega}(X)$ of idempotent probability measures with
decomposition (2) there exists an idempotent probability measure
$\xi_{12}\in\Lambda_{12}$ such that $$H(\mu_1,\
\mu_2)=\bigoplus\limits_{(x_j, x_k)\in S\xi_{12}}
|\lambda_{2k}-\lambda_{1j}|\odot\rho(x_{1j}, x_{2k}).$$}

Proof is followed from the finiteness of the set (3).

The following Lemma is rather obvious.

\textbf{Lemma 2.} \textsl{Let $\mu_1,\ \mu_2\in I_{\omega}(X)$ be an
arbitrary pair of idempotent probability measures with
decompositions (2). Then for each idempotent probability measure
$\xi_{12}\in\Lambda_{12}$ the following equations hold
$$\pi_i(S\xi_{12})=S\mu_i,\,\,\,\,\,\,\ i=1,2.$$}

Finally, we need the following technical lemma.

\textbf{Lemma 3.} \textsl{Let $\mu_1,\ \mu_2,\ \mu_3\in I_{\omega}(X)$ be
an arbitrary three of idempotent probability measures with
decompositions of the form (2). Let, moreover,
$\xi_{12}\in\Lambda_{12}$ and $\xi_{23}\in\Lambda_{23}$ be
idempotent probability measures satisfying the conclude of Lemma
1. Then there exists an idempotent probability measure
$\xi_{13}\in\Lambda_{13}$ such that $$S\xi_{13}=\{(x_{1k},x_{3l}):
\text{there is an}\ m\in\{1,...,n_2\}\ \text{such that}$$
$$(x_{1k},x_{2m})\in S\xi_{12}\ \text{and}\ (x_{2m},x_{3l})\in S\xi_{23}\}.$$}

\textsc{Proof.} Lemma 2 implies that
$\pi_2(S\xi_{12})=\pi_1(S\xi_{23})=S\mu_2$. For
$m\in\{1,...,n_2\}$ let $k_1(m),..., k_{p(m)}(m)\in\{1,...,n_1\}$
and $l_1(m),..., l_{q(m)}(m)\in\{1,...,n_3\}$ be  numbers such
that $(x_{1k_r(m)}, x_{2m})\in S\xi_{12}$ when $r=1,..., p(m)$ and
$(x_{2m}, x_{3l_s(m)})\in S\xi_{23}$ when $s=1,..., q(m)$. If is
clear that
$$\{(x_{1k_r(m)},\ x_{2m}):\ r=1,...,p(m);\ m=1,...,n_2\}=S\xi_{12},$$
$$\{(x_{2m},\ x_{3l_s(m)}):\ s=1,...,q(m);\ m=1,...,n_2\}=S\xi_{23}.$$
Consider the following functional
$$\xi=\bigoplus_{\footnotesize
\begin{array}{c}
m=1,...,n_2; \\
r=1,...,p(m); \\
s=1,...,q(m) \\
\end{array}}
(\lambda_{1k_r(m)}\odot\lambda_{3l_s(m)})\odot\delta(x_{1k_r(m)},\ x_{3l_s(m)}).
$$

Let $\lambda_{2m'}=0$ for some $m'\in\{1,...,n_2\}$. Since
$\xi_{12}$ and $\xi_{23}$ are idempotent probability measures
there exist $r\in\{1,...,p(m')\}$ and $s\in\{1,...,q(m')\}$ such
that $\xi_{1k_r(m)}=0$ and $\xi_{3l_s(m)}=0$. This means that
$\xi$ is idempotent probability measure.

If the number $m'\in\{1,...,n_2\}$ is uniquely for which the
equality $\lambda_{2m'}=0$ holds then $$\{k_1(m'),...,
k_{p(m')}(m')\}=\{1,...,n_1\}\,\ \text{and}\,\
\{l_1(m'),...,l_{s(m')}(m')\}=\{1,...,n_3\}.$$ In this case the
assertion of the Lemma carries out.

If $\lambda_{2t}=0$ fulfils only for two values $t=m'$ and $t=m''$
then $$\{k_1(m'),..., k_{p(m')}(m')\}\cup\{k_1(m''),...,
k_{p(m'')}(m'')\}=\{1,...,n_1\}$$ and $$\{l_1(m'),...,
l_{s(m')}(m')\}\cup\{l_1(m''),...,
l_{s(m'')}(m'')\}=\{1,...,n_3\}.$$ From the construction of $\xi$
follows that the assertion is true for this case too.

Continuing this process by the quantity of $t$ for which the
equality $\lambda_{2t}=0$ holds we finish the proof. Lemma 3 is
proved.

$\textbf{Proposition 1.}$ \textsl{The function
$H:I_{\omega}(X)\times I_{\omega}(X)\longrightarrow\mathbb{R}$ is
metric.}

\textsc{Proof.} Obviously that $H(\mu_1,\ \mu_2)\geq 0$ for all
pairs $\mu_1,\ \mu_2\in I_{\omega}(X)$. Let $\mu_1=\mu_2=\mu$.
Suppose $\mu$ has a decomposition
$$\mu=\lambda_1\odot\delta(x_1)\oplus
...\oplus\lambda_n\odot\delta(x_n).$$Then the idempotent
probability measures $$\xi=\lambda_1\odot\delta(x_1,x_1)\oplus
...\oplus\lambda_n\odot\delta(x_n, x_n)$$ is an element of the
$\Lambda(\mu,\ \mu)$ and for which one has $$H(\mu,\ \mu)\leq
\bigoplus\limits_{(x_j, x_j)\in S\xi}
|\lambda_{j}-\lambda_{j}|\odot\rho(x_{j}, x_{j})=0.$$ Hence
$H(\mu_1,\ \mu_2)=0$.

Inversely let now $H(\mu_1,\ \mu_2)=0$ where $\mu_1,\ \mu_2$ are
arbitrary idempotent probability measures with finite supports
admitting decompositions (2). Then (4) implies existence an
idempotent probability measures $\xi_{12}\in\Lambda_{12}$ such
that $S\xi_{12}\subset\Delta(X)\equiv\{(x,x):x\in X\}$. This is
possible only provided $S\mu_1=S\mu_2$. On the other hand again
from (4) it follows that $|\lambda_{2j}-\lambda_{1j}|=0$ for all
$j$. Hence $\mu_1=\mu_2$.

Thus for idempotent probability measures with finite supports we
have $H(\mu_1,\ \mu_2)=0$ if and only if $\mu_1=\mu_2$.

It is clear that $H$ is symmetric.

It remains to show $H$ satisfies the triangle axiom. Let $\mu_1,\
\mu_2,\ \mu_3$ are idempotent probability measures with finite
supports admitting decompositions
$$\mu_i=\lambda_{i1}\odot\delta(x_{1i})\oplus ...
\oplus\lambda_{in_i}\odot\lambda(x_{in_i}),\,\ i=1,2,3.$$ Let
$\xi_{12}$ and $\xi_{23}$ be idempotent probability measures which
exist according to Lemma 1 and $\xi'\in \Lambda(\mu_1,\ \mu_2)$ be
an idempotent probability measure existing by Lemma 3. Then
$$H(\mu_1,\ \mu_3)=\min_{
\xi\in\Lambda_{13}}\left\{\bigoplus\limits_{(x_{1k}, x_{3l})\in
S\xi} |\lambda_{3l}-\lambda_{1k}|\odot\rho(x_{1k},
x_{3l})\right\}\leq$$
$$\leq\bigoplus\limits_{(x_{1k}, x_{3l})\in S\xi'}
|\lambda_{3l}-\lambda_{1k}|\odot\rho(x_{1k}, x_{3l})\leq$$
$$\leq\bigoplus_{\footnotesize
\begin{array}{c}
(x_{1k}, x_{2m})\in S\xi_{12}, \\
(x_{2m}, x_{3l})\in S\xi_{23}  \\
\end{array}}
|\lambda_{2m}-\lambda_{1k}|\odot\rho(x_{1k}, x_{2m})\odot|\lambda_{3l}-\lambda_{2m}|\odot\rho(x_{2m}, x_{3l})\leq$$
$$\leq\bigoplus\limits_{(x_{1k}, x_{2m})\in S\xi_{12}}
|\lambda_{2m}-\lambda_{1k}|\odot\rho(x_{1k}, x_{2m})+\bigoplus\limits_{(x_{2m}, x_{3l})\in S\xi_{23}}
|\lambda_{3l}-\lambda_{2m}|\odot\rho(x_{2m}, x_{3l})=$$
$$=H(\mu_1,\ \mu_2)+H(\mu_2,\ \mu_3).$$
Proposition 1 is proved.

For idempotent probability measures $\mu_1,\ \mu_2$ with finite supports put
$$\rho_{\omega}(\mu_1,\mu_2)=\text{min}\{\text{diam}X,\ H(\mu_1, \mu_2)\}.\eqno(5)$$

Proposition 1 implies the following

\textbf{Corollary 1.} \textsl{The function
$\rho_{\omega}:I_{\omega}(X)\times
I_{\omega}(X)\longrightarrow\mathbb{R}$ is metric.}

The following proposition shows the metric $\rho_\omega$ is
expansion of the metric $\rho$.

\textbf{Proposition 2.} \textsl{Let $(X,\ \rho)$ be a compactum,
$x,\ y\in X$. Then $\rho_{\omega}(0\odot\delta_x,\
0\odot\delta_y)=\rho(x,\ y)$.}

\textsc{Proof.} Let $0\odot\delta_x$ and $0\odot\delta_y$ be
idempotent probability measures. Since $0\odot\delta_{(x,y)}\equiv
0\odot\delta_x\otimes 0\odot\delta_y$ is uniquely idempotent
probability measure lying in $\Lambda(0\odot\delta_x,\
0\odot\delta_y)$, then according to (4) and (5) one has
$$\rho_{\omega}(0\odot\delta_x,\ 0\odot\delta_y)=\rho(x,\ y).$$
Proposition 2 is proved.

Since $X$ includes into $I_{\omega}(X)$ using the natural
transformation $\delta_X:X\longrightarrow I_{\omega}(X)$, acting
as $\delta_X(x)=\delta_x\equiv 0\odot\delta_x$, then compact
Hausdorff spaces $X$ and $\delta_X(X)$ may be identify. That is
why Proposition 2 implies the following

\textbf{Corollary 2.} \textsl{For any compactum $(X, \rho)$ the
equality $\rho_{\omega}|_{X\times X}=\rho$ holds.}

Consider a sequence $\{\mu_t\}_{t=1}^{\infty}\subset
I_{\omega}(X)$ of idempotent probability measures. If the sequence
$\{\mu_t\}$ converges to some idempotent probability measure
$\mu\in I_{\omega}(X)$ concerning to metric $\rho_{\omega}$ then
it signs by $\mu_t\rightarrow \mu$. If the sequence $\{\mu_t\}$
converges to $\mu\in I_{\omega}(X)$ concerning to pointwise
convergence topology, then we use the mark $\mu_t\Rightarrow\mu$.
Recall that a sequence $\{\mu_t\}$ converges to $\mu\in
I_{\omega}(X)$ concerning to pointwise convergence topology if
provided
$\lim\limits_{t\rightarrow\infty}\mu_t(\varphi)=\mu(\varphi)$ for
all $\varphi\in C(X)$. Assume that $\mu
=\lambda_1\odot\delta_{x_1}\oplus
...\oplus\lambda_s\odot\delta_{x_s}$ and $\mu_t
=\lambda_{t1}\odot\delta_{x_{t1}}\oplus
...\oplus\lambda_{ts_t}\odot\delta_{x_{ts_t}}$, where $s$, $s_t$
are natural numbers, $t=1,2,...$.

For each positive integer $t$ consider a finite set $\{1,\ ...,\
s_t\}$. For every $t$ fix a number $i(t)\in\{1,\ ...,\ s_t\}$.
Note that then for each pair of the sequences
$\{\lambda_{ti(t)}\}_{t=1}^{\infty}\subset\mathbb{R}_{\text{max}}$
where $\bigoplus_{i=1}^{s_t}\lambda_{ti}=0$ for all $t=1,\ 2,...$,
and $\{x_{ti(t)}\}_{t=1}^{\infty}\subset X$ there exits an
idempotent probability measure $\mu_t$ such that $\mu_t
=\lambda_{t1}\odot\delta_{x_{t1}}\oplus
...\oplus\lambda_{ti(t)}\odot\delta_{x_{ti(t)}}\oplus...\oplus\lambda_{ts_t}\odot\delta_{x_{ts_t}}$.

Inversely, let $\{\mu_t\}$ be a sequence of idempotent probability
measures with finite supports. For each positive integer $t$ by an
arbitrary manner choose at one point $x_{ti(t)}\in
\text{supp}\mu_t$ and construct sequences
$\{x_{ti(t)}\}_{t=1}^{\infty}$. Similarly it may be constructed
sequences $\{\lambda_{ti(i)}\}_{t=1}^{\infty}$ of $\max$-mass of
idempotent probability measures $\mu_t$.

\textbf{Theorem 2.} \textsl{The metric $\rho_{\omega}$ generates
on $I_{\omega}(X)$ pointwise-converges topology.}

\textsc{The Proof} leans on the following two lemmas.

\textbf{Lemma 4.} \textsl{Let $\mu\in I_{\omega}(X)$,
$\{\mu_t\}_{t=1}^{\infty}\subset I_{\omega}(X)$. Then
$\mu_t\rightarrow \mu$ iff the following condition holds:}

(*) \textsl{for each point $x_i\in \text{supp}\mu$ there exists a
sequence
$\{x_{ti(t)}:x_{ti(t)}\in\text{supp}\mu_t\}_{t=1}^{\infty}$ such
that $ \rho(x_i,\ x_{ti(t)})
\mathrel{\mathop{\kern0pt\longrightarrow} \limits_{t \to \infty }}
0 $ and $ \lambda_{ti(t)}
\mathrel{\mathop{\kern0pt\longrightarrow} \limits_{t \to \infty }}
\lambda_i$ where $\lambda_{ti(t)}$ are $\text{max}$-masses of the
points $x_{ti(t)}$,  $1\leq i(t)\leq |\text{supp}\mu_t|$, $t=1,\
2,\ ...\ $.}

\textsc{Proof.} If condition (*) executes then (4) and (5)
immediately imply that $\mu_t\rightarrow\mu$.

Inversely, let $\mu_t\rightarrow\mu$. Suppose that it does not
carry out $\rho(x_i,\ x_{ti(t)})
\mathrel{\mathop{\kern0pt\longrightarrow} \limits_{t \to \infty }}
0 $, i. e. for some point $x_{i_0}\in \text{supp}\mu$ each of the
sequences
$\theta=\{x_{ti(t)}:x_{ti(t)}\in\text{supp}\mu_t\}_{t=1}^{\infty}$
does not converge to $x_{i_0}$. Then for every such sequence there
exists $\varepsilon_{\theta}>0$ such that $\rho(x_{i_0},\ x_{\tau
i(\tau)})\geq\varepsilon_\theta$ for finitely many positive
integers $\tau$.

For every $t$ by $\xi_t$ denote an idempotent probability measure
which exists by force of Lemma 1 such that $$\rho_{\omega}(\mu,
\mu_t)=\min\left\{\text{diam}\ X, \bigoplus_{(x_{i_0},\
x_{ti(t)})\in S\xi_t} |\lambda_{i_0}-\lambda_{ti(t)}|\odot
\rho(x_{i_0}, x_{ti(t)})\right\}.$$ Then for some positive
$\varepsilon$ there exists infinitely many naturals $\tau$ such
that $\rho_{\omega}(\mu, \mu_\tau)\geq \varepsilon$. This
contradicts to convergence of $\{\mu_t\}$ to $\mu$ according to
metric $\rho_{\omega}$.

Assume $\lambda_{ti(t)} \mathrel{\mathop{\kern0pt\longrightarrow}
\limits_{t \to \infty }} \lambda_i$ is false. In other words there
exists a $\max$-mass $\lambda_{i_0}$ of idempotent probability
measure $\mu$ such that each of the sequences
$\{\lambda_{ti(t)}:\lambda_{ti(t)}\ \text{is a}\ \max\mbox{-mass
of the measure}\ \mu_t\}_{t=1}^{\infty}$ does not converge to
$\lambda_{i_0}$. In this case similarly it may be taken a
contradiction with convergency of sequence $\{\mu_t\}$ to $\mu$
according to metric $\rho_{\omega}$.

Hence the convergency of a sequence $\{\mu_t\}$ to $\mu$ by metric
$\rho_{\omega}$ implies condition (*).

Thus $\rho_t\rightarrow\mu$ and (*) are equivalent. Lemma 4 is
proved.

\textbf{Lemma 5.} \textsl{Let  $\mu\in I_{\omega}(X)$,
$\{\mu_t\}_{t=1}^{\infty}\subset I_{\omega}(X)$. Then
$\mu_t\Rightarrow\mu$ iff the following condition satisfies:}

(*) \textsl{for each point $x_i\in \text{supp}\mu$ there exists a
sequence
$\{x_{ti(t)}:x_{ti(t)}\in\text{supp}\mu_t\}_{t=1}^{\infty}$ such
that $ \rho(x_i,\ x_{ti(t)})
\mathrel{\mathop{\kern0pt\longrightarrow} \limits_{t \to \infty }}
0 $ and $ \lambda_{ti(t)}
\mathrel{\mathop{\kern0pt\longrightarrow} \limits_{t \to \infty }}
\lambda_i$ where $\lambda_{ti(t)}$ are $\text{max}$-masses of the
points $x_{ti(t)}$,  $1\leq i(t)\leq |\text{supp}\mu_t|$, $t=1,\
2,\ ...\ $.}

\textsc{Proof.} Let condition (*) takes place. Let moreover for
some $i_0$,  $1\leq i_0\leq s$, we have
$\lambda_{i^0}\odot\varphi(x_{i^0})=\bigoplus\limits_{i=1}^s\lambda_i\odot\varphi(x_i)=\mu(\varphi)$,
where  $\varphi\in C(X)$ is an arbitrary function. Then there
exists  $t_0$ such that $\lambda_{ti'(t)}\odot\varphi(x_{ti'(t)})=
\bigoplus\limits_{i(t)=1}^{s_t}\lambda_{ti(t)}\odot\varphi(x_{ti(t)})=\mu_t(\varphi)$
for all  $t\geq t_0$, where $x_{ti'(t)}\rightarrow x_{i^0}$  and
$\lambda_{ti'(t)}\rightarrow\lambda_{i^0}$. Therefore for all
$t\geq t_0$  the following inequality is true
$$|\mu_t(\varphi)-\mu(\varphi)|=\left|\bigoplus\limits_{i(t)=1}^{s_t}\lambda_{ti(t)}
\odot\varphi(x_{ti(t)})-\bigoplus\limits_{i=1}^{s}\lambda_{i}\odot\varphi(x_{i})\right|=
|\lambda_{ti'(t)}\odot\varphi(x_{ti'(t)})-\lambda_{i^0}\odot\varphi(x_{i^0})|\leq$$
$$\leq|\lambda_{ti'(t)}-\lambda_{i^0}|+|\varphi(x_{ti'(t)})-\varphi(x_{i^0})|.$$

As $\varphi$ is continuous (*) implies
$|\varphi(x_{ti'(t)})-\varphi(x_{i^0})|\mathrel{\mathop{\kern0pt\longrightarrow}
\limits_{t \to \infty }}0$ and
$|\lambda_{ti'(t)}-\lambda_{i^0}|\mathrel{\mathop{\kern0pt\longrightarrow}
\limits_{t \to \infty }}0$. Then it should be
$|\mu_t(\varphi)-\mu(\varphi)|\mathrel{\mathop{\kern0pt\longrightarrow}
\limits_{t \to \infty }}0$. Since $\varphi$ is an arbitrary
function then it should be  $\mu_t\Rightarrow\mu$.

Inversely,  let  $\mu_t\Rightarrow\mu$. This means that for each
$\varphi\in C(X)$ and for an arbitrary $\varepsilon>0$ there is a
natural number  $t_{\varepsilon}$, such for all $t\geq
t_{\varepsilon}$  the following inequality it takes place
$$|\mu_t(\varphi)-\mu(\varphi)|=\left|\bigoplus\limits_{i_t=1}^{s_t}\lambda_{i_t}^t\odot
\varphi(x_{i_t})-\bigoplus\limits_{i=1}^{s}\lambda_{i}\odot\varphi(x_{i})\right|<\varepsilon .$$

Suppose that (*) is false. This means for any point
$x_{i_0}\in\text{supp}\mu$ either the sequence
$\theta=\{x_{ti(t)}:x_{ti(t)}\in\text{supp}\mu_t\}_{t=1}^{\infty}$
does not converge to  $x_{i_0}$, or corresponding sequence
$\{\lambda_{ti(t)}\}_{t=1}^{\infty}$ does not converge to
$\lambda_{i_0}$.

Assume the sequence
$\theta=\{x_{ti(t)}:x_{ti(t)}\in\text{supp}\mu_t\}_{t=1}^{\infty}$
does not converge to  $x_{i_0}$. For each such sequence there is
$\varepsilon_{\theta}>0$ such that $\rho(x_{i_0}, x_{\tau
i(\tau)})\geq\varepsilon_{\theta}$  for infinitely many naturals
$\tau$. For each such $\tau$ choose a function $\varphi_{\tau}\in
C(X)$  with $\varphi_{\tau}(x_{i_0})=\text{max}\{|\lambda_i|:i=1,
..., s\}+\varepsilon_{\theta}$  and $\varphi_{\tau}(x_{\tau
i(\tau)})=0$ for all  $i(\tau)=1,2, ..., s_{\tau}$. In carrying
out these conditions one has
$\mu(\varphi_{\tau})\geq\lambda_{i_0}\odot\varphi_{\tau}(x_{i_0})>\varepsilon_{\theta}$
and $\mu_{\tau}(\varphi_{\tau})=0$. That is why
$|\mu(\varphi_{\tau})-\mu_{\tau}(\varphi_{\tau})|=|\mu(\varphi_{\tau})|>\varepsilon_{\theta}$
for infinitely many naturals $\tau$. We obtain a contradiction
with convergence of sequence $\{\mu_t\}$ to $\mu$ by pointwise
convergence topology.

No let we assume the corresponding sequences
$\{\lambda_{ti(t)}\}_{t=1}^{\infty}$ do not converge to
$\lambda_{i_0}$. It is enough to see a case when $S\mu_t=S\mu$ for
all $t=1,2,...$ . Suppose a sequence
$\{\lambda_{ti_0(t)}\}_{t=1}^{\infty}$ does not converge to
$\lambda_{i_0}$. Then there exists $\varepsilon>0$ such
$|\lambda_{i_0}-\lambda_{\tau i_0(\tau)}|\geq\varepsilon$ for
infinite many $\tau$. For each such $\tau$ choose a function
$\varphi_{\tau}\in C(X)$ with
$\varphi_{\tau}(x_{i_0})=2|\lambda_{i_0}|$ and
$\varphi_{\tau}(x_i)=0$ for all $i=1,2,...,s$. In fulfilment these
equations one has $\mu(\varphi_{\tau})=|\lambda_{i_0}|$ and
$\mu_{\tau}(\varphi_{\tau})=0\oplus(2|\lambda_{i_0}|+\lambda_{\tau
i_0(\tau)})$. If $\mu_{\tau}(\varphi_{\tau})=0$ then
$$|\mu(\varphi_{\tau})-\mu_{\tau}(\varphi_{\tau})|=|\mu(\varphi_{\tau})|=
|\lambda_{i_0}|=|\lambda_{i_0}-\lambda_{\tau
i_0(\tau)}|\geq\varepsilon_{\theta}.$$ In the case
$\mu_{\tau}(\varphi_{\tau})=2|\lambda_{i_0}|+\lambda_{\tau
i_0(\tau)}$ it takes place the following relations
$$|\mu(\varphi_{\tau})-\mu_{\tau}(\varphi_{\tau})|=||\lambda_{i_0}|+\lambda_{\tau
i_0(\tau)}|=|\lambda_{i_0}-\lambda_{\tau
i_0(\tau)}|\geq\varepsilon_{\theta}.$$ These two inequalities
contradict the convergency of sequence $\{\mu_t\}$ to $\mu$
according to pointwise convergence topology.

Thus conditions $\mu_t\Rightarrow\mu$ and (*) are equivalent.
Lemma 5, and thereby Theorem 2 are proved.

Let $\mu\in I_{\omega}(X)$ and $\{\mu_t\}_{t=1}^{\infty}\subset
I_{\omega}(X)$. Since for any compactum $X$ the space $I(X)$ of
idempotent probability measures is also compactum then Lemmas 4
and 5 imply the following statement.

\textbf{Corollary 3.}
$\mu_t\Rightarrow\mu\Leftrightarrow\mu_t\rightarrow\mu$.

Pointwise convergence topology on $I(X)$ denote by $p$.

Since $I_{\omega}(X)$ is everywhere dense in $I(X)$ according to
$p$ for each idempotent probability measure $\mu\in I(X)\setminus
I_{\omega}(X)$ there exists $\{\mu_t\}_{t=1}^{\infty}\subset
I_{\omega}(X)$ such that $\mu_t\Rightarrow \mu$. Then
$\{\mu_t\}_{t=1}^{\infty}$ is fundamental but by force of
Corollary 3 it is nonconvergent in $I_{\omega}(X)$ according to
$\rho_{\omega}$. Thus metrical space $(I_{\omega}(X),\
\rho_{\omega})$ is not complete.

Call fundamental sequences $\{\xi_t\}_{t=1}^{\infty}$,
$\{\eta_t\}_{t=1}^{\infty}\subset I_{\omega}(X)$ are equivalent,
if $\lim\limits_{t\rightarrow\infty}\rho_{\omega}(\xi_t,
\eta_t)=0$. Consider a set $I^*(X)$ of all classes of  of
fundamental according to metric $\rho_{\omega}$ sequences in
$I_{\omega}(X)$. Let $\xi=[\{\xi_t\}_{t=1}^{\infty}]$,
$\eta=[\{\eta_t\}_{t=1}^{\infty}]\in I^*(X)$. By the equality
$$d(\xi,\eta)=\lim\limits_{t\rightarrow\infty}\rho_{\omega}(\xi_t, \eta_t)$$
we define metric on $I^*(X)$, where $\{\xi_t\}_{t=1}^{\infty}$,
$\{\eta_t\}_{t=1}^{\infty}$ are represents of corresponding
classes. It is clear that $(I^*(X), d)$ is complete metrical
space. Therefore it is unique exactly isometrical completion of
metrical space $(I_{\omega}(X), \rho_{\omega})$.

\textbf{Proposition 3.} \textsl{Spaces $(I(X), p)$ and $(I^*(X),
d)$ are homeomorphic.}

\textsc{Proof.} Let $(I^*(X), d)$ be completion of $I_{\omega}(X)$
by metric $\rho_{\omega}$. Then by density $I_{\omega}(X)$ in
$I(X)$ we have that $I(X)$ is homeomorphic put in $(I^*(X), d)$,
and since $I(X)$ is compact one has $(I(X), p)$ is homeomorphic to
$(I^*(X), d)$. Proposition 3 is proved.

Define now on $I(X)$ metric $\rho_I$ by the rule $$\rho_I(\mu,
\nu)=\lim\limits_{t\rightarrow\infty}\rho_{\omega}(\mu, \nu),\
\mu,\ \nu\in I(X), \eqno(6)$$ where $\{\mu_t\}_{t=1}^{\infty}$,
$\{\nu_t\}_{t=1}^{\infty}\subset I_{\omega}(X)$ are arbitrary
sequences such that $\mu_t\Rightarrow\mu$ and
$\nu_t\Rightarrow\nu$.

Thus we finish the proof of our main result.

Proposition 3 implies the following important statement.

\textbf{Corollary 4.} \textsl{Metric $\rho_I$ on $I(X)$ generates
pointwise convergence topology .}

Now Corollary 2 may be formulate for metric $\rho_I$.

\textbf{Corollary 5.} \textsl{$\rho_I|_{X\times X}=\rho$ for any
metrical compactum $(X, \rho)$.}

\end{document}